\documentstyle[12pt]{article}
\newcommand{\const}{\mathop{\rm const}\limits}

\newcommand{\vraisup}{\mathop{\rm vraisup}\limits}

\newcommand{\card}{\mathop{\rm card}\limits}

\newcommand{\mod}{\mathop{\rm mod}\limits}

\newcommand{\Law}{\mathop{\rm Law}\limits}

\newcommand{\Var}{\mathop{\rm Var}\limits}

\begin{document}

\begin{center}

{\bf MOMENT AND TAIL ESTIMATES FOR   MARTINGALES }\\
\vspace{3mm}
{\bf AND MARTINGALE TRANSFORM, }\\
\vspace{3mm}

{\bf with application to the martingale limit theorem in Banach spaces.}\\

\vspace{4mm}

{\bf E. Ostrovsky, L.Sirota.}\\

\vspace{3mm}

Department of Mathematics, Bar - Ilan University, Ramat - Gan, 59200,Israel.\\
e-mail: \ galo@list.ru, eugostrovsky@list.ru \\

Department of Mathematics, Bar - Ilan University, Ramat - Gan, 59200,Israel.\\
e - mail: \ sirota@zahav.net.il\\

\end{center}

{\bf Abstract.} \par

\vspace{3mm}

  In this paper non-asymptotic exponential and moment
estimates are derived for  tail of distribution for discrete time martingale
and martingale transform by means of martingale differences in the terms of {\it unconditional}
moments and tails of distributions of summands and multipliers. \par
 We show also the exactness of obtained estimations and consider some applications in
the theory of limit theorem for Banach space valued martingales. \par

\vspace{3mm}

\par

\vspace{3mm}

\vspace{4mm}
\begin{center}
\section{ Introduction. Notations. Statement of problem.} \par
\end{center}
\vspace{3mm}

 Let $ (\Omega,F,{\bf P} ) $ be a probability space,
$ \xi(1), \xi(2), \ldots,\xi(n) $ being a  centered
$ ({\bf E} \xi(i) = 0, i=1,2,\ldots,n) $ martingale - differences on the basis of the
{\it same flow }of $ \sigma - $ fields (filtration) $ F(i): F(0) =
\{\emptyset, \Omega \}, \ F(i) \subset F(i+1) \subset F, \
\xi(0) = 0; \ {\bf E} |\xi(i)| < \infty, $ and for every $ i \ge 0,
 \ \forall k = 0,1,\ldots,i -1 \ \Rightarrow  $

$$
{\bf E} \xi(i)/F(k) = 0; \ {\bf E}\xi(i)/F(i) = \xi(i) \  (\mod \ {\bf P}).
$$

\vspace{3mm}
 We denote

 $$
 S(n) = \sum_{i=1}^n \xi(i), \ n \le \infty,
 $$
we understood in the case  $ n=\infty \ S(n) $ as a limit $  S(\infty) = \lim_{n \to \infty} S(n), $
if there exists.\par
 This limit there exists if for example

 $$
 \sum_{i=1}^{\infty} \Var (\xi(i)) < \infty.
 $$

 The pair $ (S(n), F(n)) $ is (pure) martingale.\par

  Further, let $ \{ b(i) \}, i=1,2,\ldots,n $ be a {\it predictable } relatively $ \{ F(n) \} $
sequence of random variables (r.v.) such that

 $$
 \forall i \Rightarrow  {\bf E}|b(i) \xi(i)| < \infty;
 $$
then the sequence $ ( W(n), F(n) ), $ where

$$
W(n) = \sum_{i=1}^n b(i) \xi(i)
$$
is also a martingale. \par
 The transform $ S(n) \to W(n) $ is called {\it martingale transform}, generated by $ \{b(i) \}, $
or in other words, stochastic integral over discrete martingale measure. \par

\vspace{4mm}
{\bf Our aim is to obtain the moment and tail estimates for  $ S(n) $ and $ W(n) $  via the
moment and tail estimates of the sequences  $  \{ \xi(i) \} $  and $  \{ b(i) \}. $ } \par

\vspace{3mm}
 More exactly, we will estimate the distribution of $ S(n) $ and $ W(n)$ via the $ L(p) $
norms $ |\xi(i)|_p, \ |b(i)|_p $ (or via some another rearrangement invariant norms)
of a summands and multipliers  $ \xi(i), b(i), $ where we denote as ordinary for any r.v. $ \eta $

$$
|\eta|_p = \left[ {\bf E} |\eta|^p  \right]^{1/p}, \ p \in [1, \infty); \ L(p) = \{ \eta, \eta: \Omega \to R, \ |\eta|_p < \infty.  \}
$$

 Our estimates improve or generalize the well-known inequalities belonging to D.L.Burkholder \cite{Burkholder1},
 \cite{Burkholder2},  \cite{Burkholder3},  \cite{Burkholder4}, \cite{Burkholder5}; K.Bichteler \cite{Bichteler1};
J.-A.Chao  \cite{Chao1}; K.P. Choi  \cite{Choi1},  \cite{Choi2}; P. Hitczenko, S.J.Mongomery-Smith, K.Oleszkiewicz
\cite{Hitzenko1}, \cite{Hitzenko2}; A.Osekovsky \cite{Osekovsky1}, \cite{Osekovsky2}; I.Pinelis  \cite{Pinelis1},
\cite{Pinelis2}. See also the books \cite{Alon1}, \cite{Hall1}, \cite{Ledoux1}; surveys  \cite{Kallenberg1}, \cite{Peshkir1}
and articles  \cite{Astashkin1}, \cite{Astashkin2},  \cite{Azuma1}, \cite{Banuelos1}, \cite{Barral1},
\cite{Bentkus1}, \cite{Davis1}, \cite{Dzhaparidze1}, \cite{Dharmadhikari1}, \cite{Fazekas1}, \cite{Kunita1},
\cite{Laib11}, \cite{Lesign1}, \cite{Liu1}, \cite{Pena1}, \cite{Teicher1}, \cite{Zakai1} etc. \par
 Some applications of these estimates in the statistics, polymer computation, theory of percolation and theory of dynamical
systems are described in \cite{Garsia1}, \cite{Gine1}, \cite{Korolyuk1}, \cite{Laib11},
\cite{Lesign1}, \cite{Liptser1}, \cite{Yulin1}, \cite{Zhang1}. \par
 Another nearest results see in references to this work (as a rule, the last results) and in \cite{Peshkir1}.\par

\vspace{4mm}

 The paper is organized as follows. In the second section we consider a particular case when the sequence $ \{ b(i) \}  $
is non-random.  In the third section    we intend to show the exactness of our estimates up to multiplicative constant.\par

 Fourth section contains the main result of offered paper:  moments estimates for martingale transform. In the next section
we formulate and prove  some propositions about exponential tail estimate of distribution of martingale transform; we recall
before  for reader convenience some auxiliary facts about  the random variables with exponential
tails of distributions. \par

 In the sixth section  we investigate as an applications of obtained results some sufficient conditions for
weak compactness of sequence of martingale random fields, for instance, for the Central Limit Theorem in the
space of continuous functions. \par

 The last section contains some concluding remarks and generalizations. \par

 \vspace{3mm}

\section{ Moments estimates for martingales.} \par

\vspace{3mm}

{\bf Theorem 2.1.} Let $ \forall i \ \xi(i) \in L(p), \ p \ge 2. $ Then

$$
\left|n^{-1/2} S(n) \right|_p \le (p-1) \left\{ \ n^{-1} \ \sum_{i=1}^n |\xi(i)|_p^2 \right\}^{1/2}. \eqno(2.1)
$$
{\bf Proof.}  Let $ \{ b(i) \} \in B $ be in time, in this section {\it nonrandom } numerical sequence  for which

 $$
 \{ b(i) \} \in B \stackrel{def}{=} \{b = b(i): \ \sum_i b^2(i) = 1.\}
 $$

 Note that it can be assumed that $ n < \infty $ and $ p > 2 $ (the case $p=2$ is  trivial) and that
 $ \forall i \le n \ b(i) \ne 0.$ Further, the sequence $ b(i)\xi(i) $ is
also a sequence of the martingale differences relative to the source
initial filtration. \par
 We have using the main result of article \cite{Ostrovsky4}, which may be obtained in turn from the 
famous  Burkholder inequality \cite{Burkholder1}, \cite{Burkholder2}:

$$
|\sum b(i)\xi(i)|_p^p \le (p-1)^p {\bf E}\left[\sum b^2(i) \xi^2(i) \right]^{p/2}.\eqno(2.2)
$$

 Substituting into (2.2) the values $ b(i) = 1/\sqrt{n} $, we obtain what was required.\par

\vspace{3mm}

{\bf Remark 2.0.} Theorem 2.1 may be obtained also from one of the result of an article
Lesign E., Volny D. \cite{Lesign1}.

 \vspace{3mm}

 {\bf Remark 2.1.} Theorem 2.1 improved one of results of the article \cite{Ostrovsky4}, where
instead the factor $ p-1 $  in (2.1) obtained the coefficient $ p \sqrt{2}. $ \par

 \vspace{3mm}

 {\bf Remark 2.2.} It is proved in \cite{Pemantle1} that if for the martingale $ (M_n,F_n), M_1 = 0 $
 the following condition holds: $ p = \const \ge 2 \ \Rightarrow $
$$
\sup_{n \ge 2} \ \vraisup \ {\bf E} \left( |S(n) - S(n-1)|^p/F_{n-1} \right) \le Q^p < \infty, \ Q = \const < \infty
$$
then
$$
\sup_n |n^{-1/2} \ S(n)|_p \le p \ Q.
$$

\vspace{3mm}

\section{ Exactness of our estimates.} \par

\vspace{3mm}

{\bf Corollary 3.1} Denote
$$
M(p)=\sup \sup_n \sup_{b \in B} |\sum_{b \in B} b(i)\xi(i)|_p/\mu(p),
$$
where the upper bound is calculated over all the sequences of centered martingale differences
$ \{\xi(i)\} $ with finite {\it uniform} absolute moments $  \mu(p) $ of the order $p:$

$$
\mu(p) = \sup_i |\xi(i)|_p.
$$
 It follows from
(2.1) that $ M(p) \le p-1. $ On the other hands, if independent symmetrical identically distributed are considered
instead $\xi(i),$ it is proved in \cite{Ostrovsky6} that for them the fraction
in the right-hand part can have an estimate from below of the form
$ 0.87 p/\log p. $ Thus
$$
0.87 \ p/\log p \le M(p) \le  p-1, p \ge 2.
$$
 Therefore, our estimation cannot be improved essentially. \par

 Let us denote the optimal constant in the inequality (2.1) as $ K(p). $ More detail:

$$
K(p) := \sup_n \sup_{ \{\xi(i) \} }
\left[ \frac{ \left|n^{-1/2} S(n) \right|_p }{ (p-1) \left\{ \ n^{-1} \ \sum_{i=1}^n |\xi(i)|_p^2 \right\}^{1/2}} \right],
\eqno(3.1)
$$
where interior supremum in (3.1) is calculated over all centered martingale differences $ \{\xi(i) \} $ from
the space $ L(p), $ where $ p = \const \ge 2. $

{\bf Theorem 3.1.}

$$
\sup_{p \in [2,\infty]} \left[ \frac{K(p)}{p-1} \right] = 1.  \eqno(3.2)
$$
{\bf Proof.} The upper bound obtained in theorem (2.1); the lower bound in (3.2) is attained, for
instance, when $ p=2 $ and if $ \{ \xi(i) \}  $ are centered identically distributed r.v. with finite
positive variance. \par

 But the result of theorem 3.1 is not very interest,  as long as by our opinion it is very interest
to investigate the asymptotical behavior of the function $ K(p) $ as $ p \to \infty. $ \par

 Let us introduce the following constant:

$$
C = \frac{1/e}{ \left[20 \ \log^2 (2)/9 + 1/3 \right]^{1/2}} \approx 0.31080315... \eqno(3.3)
$$

{\bf Theorem 3.2}

$$
\overline{K} := \overline{\lim}_{p \to \infty} \frac{K(p)}{p-1} \ge C. \eqno(3.4)
$$
{\bf Proof.} Let us consider the following example:  $ \Omega = (0,1) $  without diadic-rational points,
$ F  \ $  is Borelian sigma field, $ {\bf P}  $ is usually Lebesgue measure.  We define a functions

$$
f(x) = |\log x| - 1; \ F(x) = \int_0^x f(t) dt = x |\log x|, \ x \in \Omega, \eqno(3.5)
$$
so that

$$
\int_{\Omega} f(x) dx = F(1-0) = 0.
$$

 We find by direct calculation using Stirling's formula as $ p \to \infty: $

$$
|S(\infty)|_p = \left[ \int_0^1 |f(x)|^p dx  \right]^{1/p} \sim \left[ \int_0^1 |\log(x)|^p dx  \right]^{1/p}=
$$

$$
\left[ \Gamma(p+1) \right]^{1/p} \sim p/e.
$$

 Further, we can and will suppose without loss of generality that the number $ p $ is integer: $ p=2,3, \ldots. $ Let us
introduce the following increasing sequence  of sigma-algebras (partitions) $ F(m), m=1,2,\ldots, $ depending on the $ p: $

$$
F(m) = \sigma \left\{\left( \frac{k}{2^{mp}},  \frac{k+1}{2^{mp}} \right)  \right\}, \
k = 0,1,\ldots, 2^{mp}-1,
$$
$ F_0 = \{\emptyset, \Omega  \}, \ F_{\infty} = F. $ \par

 We define the following (regular) martingale  $ (S(m), F(m)), \ m = 0,1,2,\ldots, \infty: $

 $$
 S(m) = {\bf E} f/F(m), \ S(0) = {\bf E} f/F(0) = {\bf E} f = 0, \ S(\infty) = f.
 $$

 We denote as usually by for any set $ A \ I(A) = I(A,x) = 1, \ x \in A, \ I(A) = I(A,x) = 0, \ x \notin A $ the indicator
function of the set $ A, $ and

$$
A_k^{(m)} =  \left( \frac{k}{2^{mp}},  \frac{k+1}{2^{mp}} \right).
$$

 The function $ S(m) $ has a view

 $$
 S(m) = \sum_{k=0}^{2^{mp}-1} 2^{mp} \cdot \ I(A_k^{(m)}, x) \cdot
 \left[ F \left( \frac{k+1}{2^{mp}} \right) - F \left( \frac{k}{2^{mp}} \right) \right],
 $$
 therefore

 $$
 S(m+1) = \sum_{l=0}^{2^{(m+1)p}-1} 2^{(m+1)p} \cdot \ I(A_l^{(m+1)}, x) \cdot
 \left[ F \left( \frac{l+1}{2^{(m+1)p}} \right) - F \left( \frac{l}{2^{(m+1)p}} \right) \right],
 $$
 and as before $ \xi(m) = S(m+1) - S(m). $   We have: $ \xi(m) = \Sigma_2(m) +  2^{mp} \times $

$$
 \sum_{l=0}^{2^{mp}-1} I(A_l^{(m+1)}, x) \cdot
 \left\{  2^{lp} \cdot \left[ F \left( \frac{l+1}{2^{(m+1)p}} \right) - F \left( \frac{l}{2^{(m+1)p}} \right) \right] -
 \left[ F \left( \frac{2l+1}{2^{mp}} \right) - F \left( \frac{2l}{2^{mp}} \right) \right] \right\}
$$

$$
=: \Sigma_2 + \sum_{l=0}^{2^{mp}-1} I(A_l^{(m+1)}, x) \cdot \Delta_1(l,m, p) \stackrel{def}{=} \Sigma_1(m) + \Sigma_2(m),
$$

$ \Sigma_2(m) =  2^{mp} \times   $

$$
 \sum_{l=0}^{2^{mp}-1} I(A_l^{(m+1)}, x) \cdot
 \left\{  2^{lp} \cdot \left[ F \left( \frac{l+2}{2^{(m+1)p}} \right) - F \left( \frac{l+1}{2^{(m+1)p}} \right) \right] -
 \left[ F \left( \frac{2l+1}{2^{mp}} \right) - F \left( \frac{2l}{2^{mp}} \right) \right] \right\}
$$

$$
= \sum_{l=0}^{2^{mp}-1} I(A_l^{(m+1)}, x) \cdot \Delta_2(l,m, p) \stackrel{def}{=} \Sigma_1(m) + \Sigma_2(m),
$$
where

$$
\Delta_1(l,m, p) = 2^{mp} \times
\left\{  2^{lp} \cdot \left[ F \left( \frac{l+1}{2^{(m+1)p}} \right) - F \left( \frac{l}{2^{(m+1)p}} \right) \right] -
 \left[ F \left( \frac{2l+1}{2^{mp}} \right) - F \left( \frac{2l}{2^{mp}} \right) \right] \right\}
$$
and analogously
$$
\Delta_2(l,m, p) = 2^{mp} \times
 \left\{  2^{lp} \cdot \left[ F \left( \frac{l+2}{2^{(m+1)p}} \right) - F \left( \frac{l+1}{2^{(m+1)p}} \right) \right] -
 \left[ F \left( \frac{2l+1}{2^{mp}} \right) - F \left( \frac{2l}{2^{mp}} \right) \right] \right\}.
$$

  We find using the explicit view of the function $ F:$

 $$
 |\Delta(l,m, p)|^2 \le 0.25 \cdot 2^{-mp} /l^p, \ l = 1,2,\ldots;
  $$

 $$
 |\Delta(0,m, p)|^2 \le 0.25 \cdot \log^2(2) \cdot m \cdot 2^{-2m}.
 $$

   The summands for $ \Sigma_2 $ are estimated analogously (moreover, they are less than ones in
 $ \Sigma_1), $  and we conclude ultimately using twice elementary inequality $ (a+b)^2 \le 2(a^2 + b^2): $
 $$
 \sum_{m=1}^{\infty} |\xi(m)|_p^2 \le
 \sum_{m=1}^{\infty} \left[\log^2(2) \cdot m^2 \cdot 2^{-2m} + 2^{-2m} \zeta^{2/p}(p) \right] =
 $$

$$
20 \ \log^2(2)/9 + (1/3) \cdot \zeta^{2/p}(p),
$$
where $ \zeta(p) $ denotes the classical Rieman zeta-function function:

$$
\zeta(p) = \sum_{k=1}^{\infty} k^{-p}, \ p > 1.
$$

 Therefore,

$$
\overline{K} \ge  \overline{\lim}_{p \to \infty} \frac{p/e}{p \cdot \left[20 \ \log^2 (2)/9 + (1/3) \cdot \zeta^{2/p}(p) \right]^{1/2}} =
$$

$$
\overline{\lim}_{p \to \infty} \frac{1/e}{ \left[20 \ \log^2 (2)/9 + (1/3) \zeta^{2/p}(p) \right]^{1/2}} =
$$

$$
\lim_{p \to \infty} \frac{1/e}{ \left[20 \ \log^2 (2)/9 + (1/3) \zeta^{2/p}(p) \right]^{1/2}} = C,
$$
as long as

$$
\lim_{p \to \infty} \zeta^{2/p}(p) = 1.
$$

\vspace{3mm}

\section{ Moments estimates for martingale transform.} \par

\vspace{3mm}

 We return in this section to the martingale transform $ S(n) \to W(n) $ estimate. Recall that
 $ \vec{b} = \{ b(i) \}, i=1,2,\ldots,n $ be here a {\it predictable} relatively $ \{ F(n) \} $  sequence of
 random variables (r.v.) such that

 $$
 \forall i \Rightarrow  {\bf E}|b(i) \xi(i)| < \infty;
 $$
then the sequence $ ( W(n), F(n) ), $ where

$$
W(n) = \sum_{i=1}^n b(i) \xi(i)
$$
is also a martingale.\par
 Let us introduce some new notations.

 $$
 |\vec{b}|_{p,\lambda} = |\vec{b}|_{p,\lambda}^{(n)} = \left[ n^{-1} \sum_{i=1}^n |b(i)|_p^{\lambda}  \right]^{1/\lambda}
 \eqno(4.1)
 $$
in the case $ n < \infty $ and

$$
 |\vec{b}|_{p,\lambda}^{(\infty)} = \sup_n |\vec{b}|_{p,\lambda}^{(n)}
$$
otherwise. Here $ p,\lambda = \const \ge 1 $ with obviously generalization when $ p=\infty $ or  $ \lambda =\infty $ or
simultaneously $ p=\infty, \lambda =\infty ; $ for instance,

$$
 |\vec{b}|_{\infty,\infty}^{(\infty)} = \sup_i \vraisup|b(i)|.
$$
 Analogously may be defined the value $ |\vec{\xi}|_{p,\mu}^{(n)}. $ \par

 {\bf Theorem 4.1.} Let $ \alpha, \beta, \lambda, \mu $ be some numbers such that

$$
( \alpha, \beta, \lambda, \mu)  \in D,
$$
 where $ D $ is the set of real number  $ D = \{ \alpha, \beta, \lambda, \mu \} $ for which
 $$
 \alpha, \beta, \lambda, \mu  \in [1, \infty], \ 1/\alpha + 1/\beta = 1, 1/\lambda + 1/\mu = 1 \eqno(4.2)
 $$
and let $ p \ge 2. $ Proposition:

$$
|W(n)|_p \le (p-1) \cdot   |\vec{b}|_{\alpha p, 2\lambda}^{(n)} \cdot |\vec{\xi}|_{\beta p,2\mu}^{(n)}. \eqno(4.3)
$$
{\bf Consequence: }

$$
|W(n)|_p \le (p-1) \cdot  \inf_{( \alpha, \beta, \lambda, \mu)  \in D  }
 \left\{ |\vec{b}|_{\alpha p, 2\lambda}^{(n)} \cdot |\vec{\xi}|_{\beta p,\mu}^{(n)} \right\}. \eqno(4.4)
$$
 {\bf Proof.} It is sufficient to consider the case $ n < \infty. $ Further, since
the sequence $ ( W(n), F(n) ) $ is also a martingale with correspondent martingale differences
$ b(i) \ \xi(i), $ we can use theorem 2.1:

$$
| n^{-1/2} \ W(n)|_p^2 \le (p-1)^2 \ n^{-1} \sum_{i=1}^n | b(i) \xi(i)|_p^2.
$$
 It follows from H\"older inequality

 $$
 | b(i) \xi(i)|_p \le |b(i)|_{\alpha p} \ |\xi(i)|_{\beta p},
 $$
following

$$
| n^{-1/2} \ W(n)|_p^2 \le (p-1)^2 \ n^{-1} \sum_{i=1}^n |b(i)|^2_{\alpha p} \ |\xi(i)|^2_{\beta p}.\eqno(4.5)
$$

 Let us introduce the following {\it normalized} measure on the {\it finite } set $ N = [1,2,\ldots, n]: $

 $$
 \nu(A) = n^{-1} \sum_{i \in A} 1 =  n^{-1} \card(A),
 $$
 then the inequality (4.5) may be rewritten as follows:

$$
| n^{-1/2} \ W(n)|_p^2 \le (p-1)^2 \  \int_N
 |b(i)|^2_{\alpha p} \ |\xi(i)|^2_{\beta p} \ \nu(di). \eqno(4.6)
$$
  The assertion of theorem 4.1  follows from (4.6) after applying H\"older inequality with powers $ \lambda, \mu. $ \par

{\bf Example 4.1.} Let the predictable sequence $ \{b(i)\} $   be bounded:

$$
V:= \sup_i \vraisup|b(i)| < \infty,
$$
then

$$
|W(n)|_p \le (p-1) \cdot V  \cdot |\vec{\xi}|_{p}^{(n)} =
(p-1) \cdot V  \cdot \left[ n^{-1} \sum_{i=1}^n |\xi(i)|_p^2  \right]^{1/2}. \eqno(4.7)
$$
 This result improved the well-known estimations belonging to D.L.Burkholder
 \cite{Burkholder1}, \cite{Burkholder2}, \cite{Burkholder3}.\par

{\bf Example 4.2.}

$$
|W(n)|_p \le (p-1) \cdot  \left[n^{-1} \sum_{i=1}^n |b(i)|^4_{2p} \right]^{1/4} \cdot
 \left[ n^{-1} \sum_{i=1}^n |\xi(i)|^4_{2p}  \right]^{1/4}. \eqno(4.8)
$$

{\bf Remark 4.1.} The estimates (4.3), (4.4) are asymptotically exact  up to multiplicative
constant still for the non-random sequence $ \{b(i) \}, $ see the third section.\par

\vspace{3mm}

\section{ Exponential tail estimate of distribution of martingale transform. }

\vspace{3mm}

 We intend to obtain in this section the exponential estimates for martingale transform,
or on the other words,  estimate of martingale transform in the Grand Lebesgue Norm.

  Let us recall a so-called "moment norm", or a norm in the Grand Lebesgue Space (GLS)  $ G(\psi) $
 on the set of r.v. defined in our probability space by the following way:
 the space $ G\psi = G(\psi) $ consist, by definition, on all the centered r.v. with finite norm

$$
||\xi||G\psi \stackrel{def}{=} \sup_{p \in [2,a)} [ |\xi|_p/\psi(p)], \ |\xi|_p =
{\bf E}^{1/p} |\xi|^p. \eqno(5.1)
$$
 Here $ \psi = \psi(p), \ p \in [2,a), \ a = \const \in (2,\infty] $ is continuous in {\it semi-open}
interval positive functions. We can in sequel conclude  that $ \forall p > a \ \psi(p) = \infty $
and $ C/\infty = 0. $  \par

 Note that the definition (5.1) is correct  still for the non-centered random
variables $ \xi.$

 Evidently, the space $ G\psi $ is rearrangement invariant (symmetrical) in the classical sense,
see for example the classical books \cite{Bennet1}, \cite{Krein1}. \par

 Recently appear many publications about these spaces, see for example
\cite{Kozatchenko1}, \cite{Fiorenza1}, \cite{Fiorenza2}, \cite{Fiorenza3},
\cite{Iwaniec1}, \cite{Iwaniec2}, \cite{Ostrovsky2}, \cite{Ostrovsky3},  \cite{Ostrovsky4},
\cite{Ostrovsky5} etc.   This spaces are convenient, e.g., for investigation of the r.v. with
exponential decreasing tail of distribution.  Indeed,
if for some non-zero r.v. $ \xi \ $ we have $ 0 < ||\xi||G(\psi) < \infty, $
 then for all positive values $ u $

$$
{\bf P}(|\xi| > u) \le 2 \ \exp \left( - \overline{\psi}^*(\log x /||\xi||G(\psi))  \right),
\eqno(5.2)
$$
where  $ \overline{\psi}(p) = p \ \log \psi(p) $ and the symbol $ g^* $ denotes some modification
of the Young-Fenchel, or Legendre  transform of the function $ g: $

$$
g^*(y) = \sup_{x \ge 2} (x y - g(x)).
$$
 see \cite{Kozatchenko1}, \cite{Ostrovsky3}, chapters 1,2.\par
 As a  consequence: if

 $$
\forall x > e^2 \ \Rightarrow  \overline{\psi}^*(\log x)  > 0,
 $$
then the space $ G\psi $ coincides with {\it exponential} Orlicz's  space over  our probabilistic space
$ (\Omega,F,{\bf P} ) $ with $ N- $ function of a view

$$
N(u) = \exp( \overline{\psi}^*(\log |u|) ), \ |u| > e^2; \ N(u) = C \cdot u^2, |u| \le e^2.
$$

  Conversely: if a r.v. $ \xi $ satisfies (5.2), then $ ||\xi||G(\psi) < \infty. $ \par

 So, the theory of $  G\psi $ spaces of random variables gives a very convenient  apparatus
for investigation of a random variables with exponential decreasing tails of distribution. \par

{\bf Remark 5.1.} If we introduce the {\it discontinuous} function

$$
\psi_r(p) = 1, \ p = r; \psi_r(p) = \infty, \ p \ne r, \ p,r \in (a,b)
$$
and define formally  $ C/\infty = 0, \ C = \const \in R^1, $ then  the norm
in the space $ G(\psi_r) $ coincides with the $ L_r $ norm:

$$
||f||G(\psi_r) = |f|_r.
$$

 Thus, the Bilateral Grand Lebesgue spaces are direct generalization of the
classical exponential Orlicz's spaces and classical Lebesgue-Riesz spaces $ L_r. $ \par

 Further, let  $ \eta(t), \ t \in T $ be a separable random process (field), $ T = \{t\} $
 is arbitrary set such that for some $ a = \const > 2 \  $ and $ p \in [2,a) $

$$
\psi_{\eta}(p) \stackrel{def}{=} \sup_{t \in T} |\eta(t)|_p < \infty.
$$
The function $ \psi_{\eta} = \psi_{\eta}(p) $ is called {\it natural function} for the
family $ \{ \eta(t) \}. $ \par
 In the term of the  natural function $ \psi_{\eta} = \psi_{\eta}(p) $ and the so-called
natural distance

$$
d_{\eta}(t,s) = ||\eta(t) - \eta(s)||G\psi_{\eta}
$$
it can be obtained under simple entropy condition the exponential estimation for the maximum
distribution  of the random field $ \eta(t) $ alike the estimate  5.2:

$$
{\bf P}(\sup_{t \in T}|\eta(t)| > u) \le 2 \ \exp \left( - \overline{\psi}^*(\log x /C_2||\xi||G(\psi))  \right).
$$

 Let us denote

 $$
 \psi_{\lambda,b}(p) = \left[n^{-1} \sum_{i=1}^n |b(i)|_p^{\lambda} \right]^{1/\lambda}
 =:|||\vec{b}|||_{p,\lambda},
 $$

 $$
 \nu_{\mu, \xi}(p) = \left[n^{-1} \sum_{i=1}^n |\xi(i)|_p^{\mu} \right]^{1/\mu}
  =:|||\vec{\xi}|||_{p,\mu},
 $$
i.e. the functions  $ \psi_{\lambda,b}(p) $ and $ \nu_{\mu, \xi}(p) $ are natural functions for the
random {\it vectors } $ \vec{b}, \ \vec{\xi}. $ \par

{\bf Remark 5.2.} The functions  $ \psi_{\lambda,b}(p) = |||\vec{b}|||_{p,\lambda} $
  and  $ \nu_{\mu, \xi}(p) = |||\vec{\xi}|||_{p,\mu} $ are so-called {\it mixed }
$  L_p(\Omega) \times L_{\lambda}(1,2,\ldots,n)   $  Lebesgue-Riesz norm for the correspondent random  vectors
$ \vec{b}, \ \vec{\xi} $  in the classical terminology of Besov-Ilin-Nikolskii, see \cite{Besov1}, chapter 11. \par
 This norms  are also rearrangement invariant for the random vectors in the sense that they dependent only on the
distribution (multidimensional)  of this vectors. \par
 More detail properties of multidimensional  rearrangement invariant spaces see in \cite{Ostrovsky1}.\par

 Recall that for two measurable
spaces $ (X, \cal{A},\mu), \ (Y, \cal{B}, \nu) $ and for two numbers $ p,q: \ 1 \le p,q \le \infty $ the mixed
$ L_{p,q} $ norm for bi-measurable real (or complex) function $ f(x,y) $ is defined as follows:

$$
|f(\cdot, \cdot)|_{p,q} = \left\{\int_Y \left[\int_X |f(x,y)|^p \ \mu(dx) \ \right]^{q/p} \ \nu(dy)  \right\}^{1/q}
$$
with obvious extension into the cases $ p = \infty $ or $ q = \infty $ or both the cases  $ p = \infty, \ q = \infty. $ \par
 This spaces are used in \cite{Besov1} in the functional analysis (imbedding theorem), in the theory of approximation etc. \par

 Further, we define a function
$$
\theta(p) \stackrel{def}{=} \inf_{(\alpha,\beta,\lambda,\mu) \in D} \left[ (p-1)\cdot \psi_{2 \lambda,b}(p) \cdot \nu_{2 \mu,\xi}(p) \right].
$$

{\bf Theorem 5.1.} Suppose  for some $ a > 2 \ \theta(a) < \infty. $ Then

$$
|| W(n) ||G\theta \le 1. \eqno(5.3)
$$
{\bf Proof} follows immediately from  (4.4) and the definition  of $ G\theta $ norm:

$$
|W(n)|_p \le (p-1) \cdot  \inf_{( \alpha, \beta, \lambda, \mu)  \in D  }
\left\{  |\vec{b}|_{\alpha p, 2\lambda}^{(n)} \cdot |\vec{\xi}|_{\beta p,2\mu}^{(n)} \right\} =
$$

$$
\inf_{(\alpha,\beta,\lambda,\mu) \in D} \left[ (p-1)\cdot \psi_{2 \lambda,b}(p) \cdot \nu_{2 \mu,\xi}(p) \right]= \theta(p).
$$

\vspace{3mm}

{\bf Example 5.1.} Assume in addition to the theorem 5.1 that all the variables $ \{b(i)\}, \ \{\xi(i) \} $ have
Gaussian distribution:

$$
\Law(b(i)) = N(0,\sigma^2(i)), \  \Law(\xi(i)) = N(0,\rho^2(i)),
$$
where
$$
0 < \min \left( \inf_i \sigma^2(i), \inf_i \rho^2(i) \right) \le  \sup_i ( \sigma^2(i) + \rho^2(i)) < \infty.
$$

 As long as

 $$
 |\xi(i)|_p \asymp \sqrt{p}, \ |b(i)|_p \asymp \sqrt{p}, \ p \in (2,\infty), \eqno(5.4)
 $$
it follows from theorem 5.1 for the values $ u \ge 1 $
$$
{\bf P}(n^{-1/2} |W(n)| > u) \le 2 \exp \left(- C_4 \sqrt{u} \right), \eqno(5.5)
$$
but really

 $$
{\bf P}(n^{-1/2} |W(n)| > u) \le 2 \exp \left(- C_5 u \right).
$$
 Notice that the equality (5.4) is true still for the uniform subgaussian random variables
$  b(i), \ \xi(i), $ i.e. for which

$$
0 < \min \left( \inf_i ( ||b(i)||G\psi_{(2)}, \inf_i ||\xi(i)||G\psi_{(2)} ) \right) \le
\sup_i ( ||b(i)||G\psi_{(2)} + ||\xi(i)||G\psi_{(2)} ) < \infty,
$$
where $ \psi_{(2)}(p) := \sqrt{p}, \ 2 \le p < \infty, $ or equally

$$
2 \exp(C_1 \lambda^2) \le {\bf E} [ \exp(\lambda b(i)) +  \exp(\lambda \xi(i)) ] \le
2 \exp(C_2 \lambda^2), \ \lambda \in R, \ C_1, C_2 = \const > 0.
$$

\vspace{3mm}

\section{ Weak compactness of sequence of martingale random fields. }

\vspace{3mm}

  We consider in this section the case when the sequence $ S(n) = S(n,v) $  dependent on some
parameter $ v; \ v \in V, \ V $ is arbitrary set. We will study the continuity and weak
compactness in the space of continuous functions $ C(V, d), $  where $ d = d(v_1,v_2) $
is some distance, the sequence of martingale random fields
$$
\overline{S}(n,v) = n^{-1/2} S(n,v)
$$
under  classical norming sequence $ 1/\sqrt{n}. $ \par

{\bf 1. Continuity.} \par
 Let $ \eta(v), v \in V $ be separable random field (r.f.) (process) defined aside
from the probabilistic space on any set $ V.$ We suppose that for arbitrary
point $ v \in V $ the r.v. $ \eta(v) $ satisfies the condition

$$
\sup_{v \in V} || \eta(v) || G\psi < \infty \eqno(6.1)
$$
for some function $ \psi = \psi(p). $  For instance, the function $  \psi(\cdot) $
may be natural function for the  field $ \eta(v), $ if there exists and is non-trivial:
$ \exists a > 2, \ \psi(a) < \infty. $ \par

 The so-called {\it natural} distance $ d(v_1, v_2) $ (more exactly, semi-distance: from
the equality $ d(v_1, v_2) = 0 $ does not follow $ v_1 = v_2 ) $ may be defined by the
formula

$$
d(v_1, v_2) = ||\eta(v_1) - \eta(v_2)||G\psi.  \eqno(6.2)
$$

The boundedness of $ d(v_1, v_2) $ follows immediately from (6.1).\par

{\bf Remark 6.1.} The continuity of the r.f. $ \eta(v) $ is understood
relative the distance $ d = d(v_1, v_2). $ \par
 We denote as usually the {\it metric entropy} of the set $ V $ in the distance $ d(\cdot, \cdot) $
as a point $ \epsilon $ as $ H(V, d, \epsilon); $ recall that $ H(V, d, \epsilon) $ is the
natural logarithm of the minimal number of $ d - $  closed balls with radius $ \epsilon, , \epsilon > 0 $
which cover the set $ V. $ By definition, $ N(V, d, \epsilon) = \exp [H(V, d, \epsilon)]. $ \par

 A very simple estimations of the values  $ N(V, d, \epsilon) $ see, e.g. in the monographs 
 \cite{Ostrovsky3}, chapter 3; \cite{Vitushkin1}. \par

 The classical theorem of Hausdorff tell us that $ \forall \epsilon > 0 \ N(V, d, \epsilon) < \infty $ iff the
set $ V $ is precompact set relative the distance $ d. $ \par
{\sc We will suppose further without loss of
generality that the set $ V $ is compact set relative the distance} $ d. $ \par

 Let us denote

$$
\psi_*(x) = \inf_{y \in (0,1)} ( x y + \log \psi(1/y) ),
$$

$$
D = \sup_{t,s \in V} d(t,s), \ H(\epsilon) = H(V,d,\epsilon).
$$

 We will use the following result \cite{Ostrovsky3}, p. 171-175:
\vspace{3mm}

{\bf Theorem 6.1.} If the following integral converges:

$$
\int_0^1 \exp(\psi_*(\log 2 + H(\epsilon)) \ d \epsilon < \infty, \eqno(6.3)
$$
then the trajectories $ \eta(v) $ are $ d(\cdot, \cdot) $ continuous with probability one:

$$
P(\eta(·) \in C(V, d)) = 1 \eqno(6.4)
$$
and moreover

$$
|| \ \sup_{v \in V}  |\eta(v)| \ ||G\psi = C_1 < \infty. \eqno(6.5)
$$

{\bf Remark 6.2.}  The case when
$$
\sup_{v \in V} |\eta(v)|_r < \infty, \ \exists r = \const \ge 1
$$
and the distance

$$
d_r(v_1, v_2) = |\eta(v_1) - \eta(v_2)|_r
$$
was considered by G.Pizier \cite{Pizier1}. Indeed, if

$$
\exists v_0 \in V, \ |\eta(v_0)|_r < \infty
$$
and
$$
\int_0^1 N^{1/r}(V, d_r, z) \ dz < \infty, \eqno(6.6)
$$
then

$$
P(\eta(\cdot) \in C(V,d_r)) = 1
$$
and

$$
| \sup_{v \in V} |\eta(v)| \ |_r < \infty. \eqno(6.7)
$$

 Notice that this result is essentially non-improved and generalized the 
classical  result belonging to A.N.Kolmogorov-Yu.V.Slutzky: if $ V = [0,1] $ and

$$
|\eta(v_1) - \eta(v_2)|_p \le C_5 \ |v_1 - v_2|^{1 + \delta}, \ \exists p \ge 1, \exists \delta > 0,
$$ 
then $ {\bf P}(\eta(\cdot) \in C[0,1]) = 1. $\par

\vspace{3mm}

 Let now $ \eta_n(v), v \in V, \ n=1,2,\ldots $ be a
{\it family} of separable random fields (r.f.) (processes) defined aside
from the probabilistic space on any set $ V.$ We suppose that for arbitrary
point $ v \in V $ the r.v. $ \eta_n(v) $ satisfies the condition

$$
\sup_n \sup_{v \in V} || \eta_n(v) || G\psi < \infty \eqno(6.8)
$$
for some function $ \psi = \psi(p). $  For instance, the function $  \psi(\cdot) $
may be natural function for the random fields $ \eta_n(v), $ if there exists and is non-trivial:
$ \exists a > 2, \ \nu(a) < \infty, $ where

$$
\nu(p) = \sup_n \sup_{v \in V} |\eta_n(v)|_p.
$$

 The so-called natural distance $ d_{\infty}(v_1, v_2) $ (more exactly, semi-distance) in
the considered   case of the family of separable random fields (r.f.) (processes) $ \eta_n(v) $
 may be defined by the formula

$$
 d_{\infty}(v_1, v_2) = \sup_n||\eta_n(v_1) - \eta_n(v_2)||G\psi.  \eqno(6.9)
$$

The boundedness of $  d_{\infty}(v_1, v_2) $ it follows immediately from (6.8).\par

{\bf Remark 6.2.} The continuity of the random fields $ \eta_n(v) $ is understood
relative the distance $ d_{\infty} = d_{\infty}(v_1, v_2). $ \par

 Let us denote as before

$$
\nu_*(x) = \inf_{y \in (0,1)} ( x y + \log \nu(1/y) ),
$$

 We will use again the following result \cite{Ostrovsky3}, p. 171-175:
\vspace{3mm}

{\bf Theorem 6.2.} If the following integral converges:

$$
\int_0^1 \exp(\nu_*(\log 2 + H(V, d_{\infty}, \epsilon)) \ d \epsilon < \infty, \eqno(6.10)
$$
and the family of distributions on the real line of one-dimensional r.v. $ \eta_n(v_0) $
for some $ v_0 = \const \in V $ is weakly compact,
then all the trajectories $ \eta_n(v), \ n=1,2,\ldots  $ are $ d_{\infty}(\cdot, \cdot) $ continuous with 
probability one:

$$
{\bf P}(\eta_n(\cdot) \in C(V, d_{\infty})) = 1 \eqno(6.11)
$$
and moreover the family of distributions $ \mu_n(\cdot) $ on the space $ C(V, d_{\infty}) $ generated by
random fields $ \eta_n(v):  $

$$
\mu_n(A)  = {\bf P} (\eta_n(\cdot) \in A), \ A \subset C(V, d_{\infty})
$$
is weakly compact. \par

\vspace{3mm}

{\bf Remark 6.3.} In the case when
$$
\sup_n \sup_{v \in V} |\eta_n(v)|_r < \infty, \ \exists r = \const \ge 1
$$
and the  correspondent distance $ d_{\infty,r}(v_1, v_2) $ is introduced  by the following way:

$$
d_{\infty,r}(v_1, v_2) = \sup_n |\eta_n(v_1) - \eta_n(v_2)|_r,
$$
then the condition (6.10) has a view:

$$
\int_0^1 N^{1/r}(V, d_{\infty,r}, z) \ dz < \infty,
$$
i.e. coincides with above mentioned Pizier's condition \cite{Pizier1}.\par

\vspace{3mm}

{\bf 2.  Weak compactness of sequence of martingale random fields.}

\vspace{3mm}

 Let us return in this pilcrow to the {\it martingale case,}  but we suppose here that the centered martingale
differences  $ \xi(i) $  dependent in addition on some parameter $ v; \ v \in V, $ where $  V $ is arbitrary set:

$$
\xi = \xi(i,v), \ v \in V.
$$

 We denote as in the second section

 $$
 \overline{S}(n,v) =  n^{-1/2} S(n,v) = n^{-1/2} \sum_{i=1}^n \xi(i,v);
 $$

$$
 \tau(p) = \sup_n  \sup_{v \in V} \left[  (p-1) \left\{ \ n^{-1} \ \sum_{i=1}^n |\xi(i,v)|_p^2 \right\}^{1/2} \right], \eqno(6.12)
$$
and suppose $ a:= \sup_{(p: \tau(p) < \infty)} p > 2 $  (may be, $ a = \infty;) $

$$
\rho(v_1,v_2) = \sup_{p \in (2,a)} \sup_{v \in V}
\left\{  \left[ (p-1) \left\{ \ n^{-1} \ \sum_{i=1}^n |\xi(i,v_1) - \xi(i,v_2)|_p^2 \right\}^{1/2} \right]/\tau(p) \right\}. \eqno(6.13)
$$

\vspace{3mm}

{\bf Theorem 6.3.} If the following integral converges:

$$
\int_0^1 \exp(\tau_*(\log 2 + H(V, \rho, \epsilon)) \ d \epsilon < \infty, \eqno(6.14)
$$
and the family of distributions on the {\it real line} of one-dimensional r.v. $ \overline{S}(n,v_0) $
for some $ v_0 = \const \in V $ is weakly compact,
then all the trajectories $ \overline{S}(n,v), \ n=1,2,\ldots  $ are $ \rho(\cdot, \cdot) $ continuous with probability one:

$$
{\bf P}( \overline{S}(n,\cdot) \in C(V, \rho))) = 1 \eqno(6.15)
$$
and moreover the family of distributions on the space $ C(V, \rho) $ generated by
random fields $ \overline{S}(n,\cdot)  $ is weakly compact. \par

\vspace{3mm}

{\bf Proof.}  It follows from theorem 2.1 that

$$
\sup_n\sup_{v \in V} \left|n^{-1/2} S(n,v) \right|_p \le
\sup_n\sup_{v \in V} \left[(p-1) \left\{ \ n^{-1} \ \sum_{i=1}^n |\xi(i,v)|_p^2 \right\}^{1/2} \right]=
\tau(p). \eqno(6.16)
$$
 Applying again the proposition of theorem 2.1 to the sequence of martingale differences
 $  \xi(i,v_1) - \xi(i,v_2), $ we obtain analogously

$$
\sup_n  \sup_{p \in (2,a)} \left| n^{-1/2} (S(n,v_1) - S(n,v_2)) \right|_p/\tau(p) \le
$$

$$
\sup_n  \sup_{p \in (2,a)}
 \left\{\left[(p-1) \left\{ \ n^{-1} \ \sum_{i=1}^n |\xi(i,v_1) - \xi(i,v_2)|_p^2 \right\}^{1/2} \right]/\tau(p)\right\}
  = \rho(v_1, v_2). \eqno(6.17)
$$

 It remains to use the proposition of theorem 6.2. \par

\vspace{3mm}

{\bf Example 6.1.} Suppose in addition of theorem 6.3

$$
\tau(p) = \psi_r(p), \ \exists r \ge 2.
$$
Then the condition (6.14) of theorem 6.3 has a view

$$
\int_0^1 N^{1/r}(V, \rho, z) \ dz < \infty. \eqno(6.18)
$$

 Further, let for instance $  V $ be bounded closed subset of the set $ R^d $ equipped 
with ordinary Euclidean distance $ |v_1 - v_2|. $   Assume that

$$
\rho(v_1,v_2) \asymp C \ |v_1 - v_2|^{\alpha}, \ \alpha = \const \in (0,1].
$$

 The condition  (6.18) is satisfied iff $ r > d/\alpha. $ \par

\vspace{3mm}

{\bf Example 6.2.} Suppose in addition to the conditions of theorem 6.3

$$
\tau(p) = \psi_{(2)}(p) = \sqrt{p}, \ 2 \le p \le \infty.
$$
Then the condition (6.14) of theorem 6.3 has in this case a view

$$
\int_0^1 H^{1/2}(V, \rho, z) \ dz < \infty. \eqno(6.19)
$$

 The condition  (6.19) is satisfied if for example

$$
 H(V, \rho, \epsilon) \ dz <  C \ \epsilon^{- 2 + \delta }, \ \epsilon \in (0,1), \ \exists \delta = \const > 0.
$$
 Notice that the condition (6.19) coincides with the famous {\it sufficient} condition belonging to
R.M.Dudley  \cite{Dudley1} and X.Fernique \cite{Fernique1} for continuity of Gaussian processes.   \par
 Since this condition (6.19) is also {\it necessary} for stationary process $ \overline{S}(\infty,v), \
 v \in [0, 2 \pi), $ we conclude that this condition is essentially non-improvable for martingale
 limit theorem in the considered case. \par

\vspace{3mm}

\section{ Concluding remarks and applications.} \par

\vspace{3mm}

{\bf A. Limit theorem for martingales in Banach space. }\par

\vspace{3mm}

 Assume in addition to the conditions of theorem 6.3 that the finite-dimensional distributions of the
sequence $  \overline{S}(n,v), \ n=1,2,\ldots  $  converge to the finite-dimensional distributions of
some non-trivial random field $ \overline{S}(\infty,v). $ The sufficient conditions for this
convergence may be find in the classical book of Hall P. and Heyde C.C. \cite{Hall1}. \par
 As a rule, the limiting field has a Gaussian distribution or multiple stochastic integral over
Gaussian stochastic measure with independent values on the disjoint sets (Non-Central
Limit Theorems.) \par
 Then the sequence of random fields   $  \overline{S}(n,v), \ n=1,2,\ldots  $ converges weakly
in the space of $ \rho- $  continuous functions $ C(V,\rho) $  to the random field
$  \overline{S}(\infty,v).  $ \par
 As a consequence: for arbitrary continuous functional $ Z, \ Z: C(V,\rho) \to R, $ for example,
$ Z(f) = \sup_{v \in V} |f(v)|  $

$$
\lim_{n \to \infty} {\bf P} (Z( \overline{S}(n,\cdot) > u)) = {\bf P} (Z( \overline{S}(\infty,\cdot) > u)).
$$
 The last equality was used, e.g., in the method Monte-Carlo  \cite{Frolov1}.\par
 In the case when $  \overline{S}(\infty,v)  $ is Gaussian, we obtain the sufficient conditions
for martingale Banach space valued Central Limit Theorem. \par

\vspace{3mm}

{\bf B. Independent case. }\par

\vspace{3mm}

 If in addition to the conditions of theorem 6.3 the r.f. $ \xi(i_1,v_1), \xi(i_2,v_2), \ldots, \xi(i_m,v_m), \ m=1,2,\ldots $
are completely independent for $ i_s \ne i_l, $ we obtain from the theorem 6.3 the classical CLT in the
space of continuous functions, see, e.g. \cite{Kozatchenko1}, \cite{Dudley2}, \cite{Ledoux1}. About the CLT
in another separable Banach spaces see  \cite{Vakhania1}. \par

\vspace{3mm}

{\bf C. The case of stochastic integrals instead sums. }\par

\vspace{3mm}
 It may be considered analogously to the section 3 the case  of stochastic integrals over continuous martingale
instead sums

$$
W(t) = \int_{(0,t)} b(s) \ d M(s), \eqno(7.1)
$$
where $ M(t) $ is left continuous square integrable martingale or semimartingale and $  b(t) $ is predictable random
process, see \cite{Burkholder2}, \cite{Burkholder5}, \cite{Osekovsky2}, \cite{Osekovsky1}. \par

\vspace{3mm}

{\bf D. Multiple martingale transform. }\par

\vspace{3mm}

 The moment and tail estimates for the  {\it multiple} martingale transform, i.e. the transform
of a view

$$
Q(d,n, \{ \xi(\cdot, \cdot, \ldots,\cdot)  \}) =
\sum_{ \ {\bf i } \  \in I(d,n)} b( \ {\bf i } \ ) \ \xi( \ {\bf i } \ ), \eqno(7.2)
$$
where

$$
{\vec {\bf i}} \in I \ \Rightarrow \xi( \ {\vec {\bf i}} \ )
\stackrel{def}{= } \prod_{s=1}^d \xi(i_s,s), \eqno(7.3)
$$

$ I = I(n) = I(d,n) = \{ i_1,i_2,\ldots,i_d \}, \ $ {\it is the set} of
indices of the form $ I(n) = I(d,n) =  \{ {\vec {\bf i}} \} = \{ {\bf i } \}  =
\{ i_1,i_2,\ldots,i_d \}, $ such that
$ 1 \le i_1 < i_2 \ldots < i_{d-1} < i_d \le n, \ $
with non-random multiple sequence $  b( \ {\bf i } \ ) $ are obtained, e.g. in
\cite{Ostrovsky4}.\par

\vspace{3mm}

{\bf D. Quadratic  $p-$ characteristic version of our inequality. }\par

In the theory of martingales the quantity

$$
 [f ]_n = \sum_{i=1}^n \xi^2(i) \eqno(7.4)
$$
is widely called the quadratic variation of the  martingale $ (S(n), F(n)),  $  see \cite{Peshkir1}.
 The classical Burkholder inequality connected the $ L(p) $ estimates between $ S(n) $ and $ [f]_n. $ \par
We introduced the new parameter, say $ p- $ quadratic variation  $ [f ]_{n,p}   $
of the  martingale $ (S(n), F(n)):  $

$$
[f ]_{n,p} = \left\{ \sum_{i=1}^n |\xi^2(i)|_p \right\}^{1/2} \eqno(7.5)
$$
and investigated  the $ L(p) $ relations  between $ S(n) $ and $ [f]_{n,p}. $ \par

 In the square-integrable case $ {\bf E} S^2(n) < \infty  $  a significant role is
played the so-called $ p- $ {\it quadratic characteristic:} $ <f >_{n,p}, $ at last in the case $ p=2: $

$$
<f>_{n,p} = \left\{ \sum_{i=1}^n   {\bf E} \left| \ \xi^2(i) /F(i-1) \ \right|_p \right\}^{1/2}. \eqno(7.6)
$$

 In this terms an analog of theorem 2.1 may formulated as follows:

$$
\left|n^{-1/2} S(n) \right|_p \le  \tilde{K}(p) <f>_{n,p}, \eqno(7.7)
$$
where for the optimal value of variable  $ \tilde{K}(p), $ namely

$$
\tilde{K}(p):= \sup_n \sup_{\{\xi(i) \} } \left[ \frac{\left|n^{-1/2} S(n) \right|_p}{ <f>_{n,p}} \right] \eqno(7.8)
$$
where the upper bound is calculated over all the sequences of centered martingale differences
$ \{\xi(i)\} $ with finite absolute moments of the order $ p, $  is valid the following double
inequality:

$$
 C_1 \ \frac{p}{\log p} \le   \tilde{K}(p) \le C_2 \ \frac{p}{\log p}, \ p \ge 2, \eqno(7.8)
$$
where $ C_1, C_2 $ are finite positive {\it absolute} constants.\par

\vspace{3mm}

 {\bf Proof } is very simple. The upper bound in the last inequality may be obtained
analogously the proof of theorem 2.1  by means of inequality  A.Osekovsky \cite{Osekovsky1}
instead Burkholder's inequality, namely:

$$
\left|n^{-1/2} S(n) \right|_p \le  C_3 \cdot [p /\log p] \cdot \left| \ <f>_{n,2} \ \right|_p, \ p \ge 2, \eqno(7.9)
$$
the lower bound is attained, for instance, for the
sequences of centered independent random variables $ \{ \xi(i) \} $ with finite absolute
moments of the order $ p, \ p \ge 2, $  see many articles \cite{Rosenthal1}, \cite{Ibragimov1},
\cite{Ibragimov2},  \cite{Hitzenko2}, \cite{Ostrovsky6}, \cite{Utev1} etc.\par
 We remain to reader to generalize the last inequality on the Grand Lebesgue Spaces.\par

\end{document}